\documentclass[letterpaper, 10 pt, conference]{ieeeconf}  % Comment this line out if you need a4paper

\IEEEoverridecommandlockouts                              % This command is only needed if
\usepackage{graphicx}
\usepackage{amsfonts,amsmath,amssymb,bm}
\usepackage{color,cases,float}
\usepackage{hyperref}
\usepackage{psfrag}
\usepackage{eucal}
\usepackage{algorithm}
\usepackage{overpic}
\usepackage{subfigure}
\usepackage{url}
\usepackage[all,cmtip]{xy}
\usepackage[noend]{algpseudocode}

\def\R{\mathbb{R}}

\newtheorem{proposition}{Proposition}
\newtheorem{lemma}{Lemma}
\newtheorem{theorem}{Theorem}
\newtheorem{corollary}{Corollary}
\newtheorem{remark}{Remark}
\newtheorem{assumption}{Assumption}

\newenvironment{proofof}[1]
  {\par\noindent\hspace{2em}{\itshape Proof of #1: }}
  {\hspace*{\fill}~\QED\par}

\title{\LARGE \bf
Stability Filters for Receding-Horizon Games
}

\author{Tatiana Tatarenko and Alberto Bemporad% <-this % stops a space
\thanks{T.~Tatarenko is with the Department of Control Theory and
Intelligent Systems, TU Darmstadt, Germany. E-mail:
{\tt\small tatiana.tatarenko@tu-darmstadt.de}.}%
\thanks{A.~Bemporad is with the IMT School for Advanced Studies,
Piazza San Francesco 19, Lucca, Italy. E-mail:
{\tt\small alberto.bemporad@imtlucca.it}.}%
\thanks{This work was funded by the Deutsche
Forschungsgemeinschaft (DFG, German Research Foundation), project
number 528033031. This work was also funded by the European Union (ERC Advanced
Research Grant COMPACT, No.~101141351). Views and opinions expressed
are however those of the authors only and do not necessarily reflect
those of the European Union or the European Research Council. Neither
the European Union nor the granting authority can be held responsible
for them. }%
}

\begin{document}

\maketitle
\thispagestyle{empty}
\pagestyle{empty}

%====================================================================
\begin{abstract}
In receding-horizon games, self-interested players repeatedly solve a
finite-horizon generalized Nash equilibrium problem (GNEP) and apply its
first-stage actions. Ensuring stability of the resulting
closed-loop system, however, is a challenging task that requires
careful design of the players' strategies or additional constraints.
We propose a stability filter that minimally modifies the variational 
equilibrium of the finite-horizon GNEP to enforce a prescribed
decrease of a quadratic Lyapunov function. For receding-horizon games 
with linear dynamics, quadratic costs, and polyhedral coupling input
constraints, whose associated GNEP is strongly monotone, 
we show that a filter always exists such 
that a closed-loop equilibrium exists, the filter is
always feasible, does not intervene when the
equilibrium already satisfies the decrease condition, and that 
renders the closed loop globally exponentially stable. Moreover, the
filtered strategy is an $\epsilon$-generalized Nash equilibrium, with
$\epsilon$ vanishing exponentially along closed-loop trajectories.
\end{abstract}

%====================================================================
\section{Introduction}
\label{sec:intro}
%====================================================================

Many control problems involve several decision makers that are acting on the same
dynamical system while pursuing individual objectives, for instance in
demand-side management \cite{chen2014autonomous} or the charging of
electric vehicles \cite{ma2011decentralized}. Coupled through the shared
dynamics and the shared constraints, their interaction is naturally modeled as a generalized
Nash equilibrium problem (GNEP)
\cite{facchinei2007generalized,fischer2014generalized}, whose
variational equilibrium (v-GNE) is characterized by a variational
inequality with equal multipliers for all players
\cite{KULKARNI2012,FaccPang1}. Receding-horizon games (RHGs) combine
this model with the receding-horizon principle of model predictive
control (MPC) \cite{VANDENBROEK2002937,hall2024receding,benenati2026lq}: at every time step, the players
choose a sequence of future control actions by solving a finite-horizon GNEP,
parameterized by the current state, and apply the first-stage actions of the 
resulting equilibrium. RHGs thus describe
strategic agents that repeatedly re-plan their strategies
in response to the evolving state.

Unlike in standard MPC, where the stage cost is chosen by the designer
to optimize all control actions simultaneously, the
costs in a RHG belong to the players and reflect their individual
interests. Hence, the standard argument of using a single {\it optimal cost} 
as a Lyapunov function fails beyond potential games. Certifying closed-loop stability
is of paramount importance for the applicability of RHGs in practice as, for example, 
agents can destabilize each other through an RHG even when their local dynamics
are stable \cite{hall2026stability}. 

For strongly monotone games with stable linear dynamics, \cite{hall2026stability} derives, via
dissipativity and monotone operator theory, a sufficient condition in
the form of a linear matrix inequality (LMI) that certifies global asymptotic
stability, together with guidelines for tuning the cost weights.
Other works equip the finite-horizon game with terminal ingredients:
\cite{benenati2026lq} obtains an explicit stabilizing terminal cost
for linear-quadratic games from an equivalent regulator in an
augmented space, and \cite{hante2024quasi,hante2025stabilizing} derive
conditions on the players' objectives that guarantee asymptotic
stability under terminal constraints or Lyapunov terminal costs via
$\alpha$-quasi GNEPs.
 %Without terminal ingredients,
%\cite{hall2026nonlinear} establishes practical asymptotic convergence
%of nonlinear RHGs from the turnpike phenomenon, and
%\cite{yildirim2026misspecification} gives stability criteria when
%agents misspecify each other's objectives.
 Closed-loop stability under input
and state constraints is instead achieved in~\cite{DFB26} by introducing in the GNEP 
a shared linear constraint on the first control moves that forces a certain polyhedral
control Lyapunov function to decrease. These approaches either certify stability of a given RHG, with no
remedy when the certificate fails, or modify the players' costs or
constraints, which may severely perturb the players' private
interests. 

A different mechanism proposed in the MPC literature to ensure stability is the use of stability filters
\cite{wabersich2021predictive,wabersich2023datadriven,hsu2024safety}
and control barrier and control Lyapunov function based quadratic
programs \cite{ames2017cbf,ames2014rapidly}, which minimally modify a
nominal input to enforce safety or stability.

\emph{Contributions.} For open-loop Schur stable linear dynamical systems, quadratic
costs defining a game with strongly monotone pseudo-gradient, and compact polyhedral
coupling constraints, we propose a stability filter that
replaces the v-GNE at every time step by the closest feasible joint
strategy that guarantees a prescribed decrease of a quadratic
Lyapunov function. Since the filtered strategy is in general no longer
an equilibrium, a central question is how much strategic optimality
is sacrificed for stability. We establish the following properties:
\begin{enumerate}
    \item The proposed convex filter is feasible for every state,
    leaves the v-GNE unchanged whenever it already satisfies the
    decrease condition, and intervenes by an amount bounded linearly
    in the distance to the equilibrium.
    \item The filtered closed loop is globally exponentially stable,
    with a rate prescribed by the designer.
    \item The filtered strategy is an $\epsilon$-generalized Nash
    equilibrium, with $\epsilon$ bounded by the v-GNE multipliers and
    the size of the intervention and vanishing exponentially along the
    closed loop. Hence strategic optimality is recovered
    asymptotically.
\end{enumerate}

\emph{Organization.} Section~\ref{sec:problem} introduces the RHG and
the standing assumptions. Section~\ref{sec:stability_filter} establishes
existence of a closed-loop equilibrium, introduces the stability
filter, proves global exponential stability of the filtered closed
loop, and analyzes the resulting strategic loss; the main results are
summarized in Theorem~\ref{thm:main}. A numerical example is provided in Section~\ref{sec:example}
and the paper concludes with a discussion in Section~\ref{sec:conclusion}.

\section{Problem Formulation}
\label{sec:problem}

We consider the linear system
\begin{equation}
    x(t+1)=Ax(t)+Bv(t),
    \qquad t=0,1,2,\dots,
    \label{eq:system}
\end{equation}
where $x(t)\in\mathbb{R}^{n_x}$ is the system state at time $t$ and
$v(t)\in\mathbb{R}^{n_v}$ is the command input
applied by all $N$ players. In game-theoretic MPC, 
$v(t)$ is obtained from the first stage of the joint finite-horizon strategy, i.e.,
$v(t)=\Xi u(t)$, where $u(t)\in\mathbb{R}^{n_u}$ denotes the joint finite-horizon strategy
generated at time $t$ by solving the receding-horizon game, and the matrix $\Xi$
extracts the first-stage joint action from the finite-horizon
strategy.

The strategy $u(t)=(u_i(t))_{i\in\mathcal N}$ is partitioned into the
decision blocks indexed by $\mathcal N=\{1,\dots,N\}$. 
Each player's block is itself a
finite-horizon sequence of length $N_p$,
\begin{equation}
    u_i(t)=(u_i(0|t),\dots,u_i(N_p-1|t)),
    \label{eq:ui_seq}
\end{equation}
where $u_i(k|t)$ is player $i$'s action planned at time $t$ for
prediction step $k$, and $u_{-i}(t)=(u_j(t))_{j\neq i}$ collects the
sequences of all other players. The cost of player $i$ is
% \begin{equation}
% \begin{split}
%     J_i(u_i,u_{-i},x)
%     =\;&
%     \sum_{k=0}^{N_p-1}
%     \Bigl[
%         \tfrac12 u_i(k|t)^\top W_{ii}u_i(k|t)
%     \Bigr.
%     \\
%     &
%         +
%         u_i(k|t)^\top\sum_{j\neq i}W_{ij}u_j(k|t)
%     \\
%     &\Bigl.
%         +
%         u_i(k|t)^\top K_i x(k|t)
%         +
%         u_i(k|t)^\top w_i
%     \Bigr],
% \end{split}
% \label{eq:Ji_quad}
% \end{equation}
\begin{equation}
    \begin{split}
    J_i(u_i,u_{-i},x) = &\sum_{k=0}^{N_p-1}
    \begin{bmatrix}u_i(k|t) \\ u_{-i}(k|t) \\ x(k|t)\end{bmatrix}^\top\!\!
    {\mathcal W}_i
    \begin{bmatrix}u_i(k|t) \\ u_{-i}(k|t) \\ x(k|t)\end{bmatrix}\\
    {\mathcal W}_i = &\begin{bmatrix}
        W_{ii} & W_{i,{-i}} & K_i \\
        W_{-i,i} & W_{-i,{-i}} & K_{-i} \\
        K_i^\top & K_{-i}^\top & Q_i
    \end{bmatrix}
    \end{split}
    \label{eq:Ji_quad}
\end{equation}
% \AB{I've removed $ u_i(k|t)^\top w_i$ from the cost function, it is an unusual term that biases the input away from zero and, therefore, the equilibrium state too. It may re-introduced by assuming that deviations of $u$ and $x$ from a steady-state reference, see~\eqref{eq:tracking}. For the ECC paper we can just focus on regulation to the origin only $r=0$.
% }
with ${\mathcal W}_i\succeq 0$ and $W_{ii}=W_{ii}^\top\succ0$
%\AB{We want $J_i(u_i,u_{-i},x)$ strongly convex wrt $u_i$ after eliminating the predicted states}, 
where the predicted states are
generated from the current state $x=x(t)=x(0|t)$ as
\begin{equation}
    x(k|t)=A^kx + \sum_{j=0}^{k-1}A^{k-1-j}Bu(j|t),
    \label{eq:prediction_model}
\end{equation}
with $u(j|t)=(u_i(j|t))_{i\in\mathcal N}$. Since the game solved at
time $t$ depends on $t$ only through $x=x(t)$, we omit the time index
in $J_i$ and in all quantities derived from it below, which we treat as functions of a generic state $x$.

Given $x$ and the sequences $u_{-i}$, each player $i$ solves
\begin{equation}
    u_i
    \in
    \operatorname*{arg\,min}_{z_i:\,(z_i,u_{-i})\in Z}
    J_i(z_i,u_{-i},x)
    \qquad
    \forall i\in\mathcal N,
    \label{eq:Ji}
\end{equation}
where the polyhedral constraint set
\begin{equation}
    Z
    =
    \left\{
        u\in\mathbb{R}^{n_u}
        \;\middle|\;
        Cu\le d,\quad Eu=e
    \right\}
    \label{eq:Z}
\end{equation}
encodes shared linear inequality and, possibly, equality constraints the
input variables of all players need to satisfy.
% \AB{Can we extend to state-dependent constraints as well, so to be able to 
% handle state constraints?
% \[
%     Z(x)=\left\{u\in\mathbb{R}^{n_u} \;\middle|\; C_uu+C_xx\le d,\quad E_uu+E_xx=e\right\}
% \]
% }

The players' pseudo-gradient $F(u,x)=(F_i(u,x))_{i\in\mathcal N}$
stacks the partial gradients
$F_i(u,x)=\nabla_{u_i}J_i(u_i,u_{-i},x)$. As from
\eqref{eq:prediction_model} $x(k|t)$ is affine in
$(u,x)$, $J_i$'s are quadratic and the associated 
pseudogradient $F$ is affine, i.e.
\begin{equation}
    F(u,x)=Hu+Gx+h,
    \label{eq:affine_F}
\end{equation}
with $H$, $G$, $h$ obtained from the standard condensed 
representation of~\eqref{eq:Ji}, and $H_i$, $G_i$, $h_i$
denoting the corresponding $i$-th block row/entry.
Note that, for better generality, we retain the constant term $h$ in the affine representation of the pseudogradient, even though $h_i=0$ for each player $i$ when the game is generated from~\eqref{eq:Ji_quad}--\eqref{eq:Z}.

We impose the following standing assumptions.
\begin{assumption}
\label{ass:main}
The following conditions hold:
\begin{enumerate}
    \item $Z$ is a nonempty and compact polyhedron;
    \item $F(\cdot,x)$ is strongly monotone on $Z$, uniformly in $x$,
          with constant $\mu_F>0$, i.e., 
$\big\langle F(z_1,x)-F(z_2,x),\, z_1-z_2\big\rangle
\ge
\mu_F\|z_1-z_2\|^2$
for all $z_1,z_2\in Z$ and all admissible $x$.
    \item for every player $i$, the cost
          $J_i(\cdot,u_{-i},x)$ is $m_i$-strongly convex in the
          player's own decision variable, with $m_i>0$;
    \item $A$ is Schur stable.
\end{enumerate}
\end{assumption}

\begin{remark}[Relation between Assumption~\ref{ass:main}.2 and~\ref{ass:main}.3]
\label{rem:mi_redundancy}
Since $F_i(\cdot,x)=\nabla_{u_i}J_i(\cdot,u_{-i},x)$, 
the strong monotonicity inequality of Assumption~\ref{ass:main}.2 
implies Assumption~\ref{ass:main}.3 with $m_i\ge\mu_F$,
as can be seen by testing the inequality on
vectors $z_1,z_2\in Z$ such that $z_{1,-i}=z_{2,-i}=0$.
We retain it as a
separate hypothesis only to allow, for each player, a constant $m_i$
tighter than $\mu_F$. Such tighter constants directly sharpen the
bounds of Theorems~\ref{thm:strategic_loss} and~\ref{thm:epsilon_decay} below. % through $c_H=\max_i\|H_i\|^2/(2m_i)$.
\end{remark}

Under Assumption~\ref{ass:main}, for a given state $x$, the
variational generalized Nash equilibrium (v-GNE) of the game
\eqref{eq:Ji} is characterized by the variational inequality
$VI(Z,F(\cdot,x)),$ which, by strong monotonicity of $F$ and Assumption~\ref{ass:main}.1, 
admits a unique solution for every $x$. We denote the corresponding equilibrium solution mapping
by
\begin{equation}
    S(x)=\operatorname{SOL}(F(\cdot,x),Z).
\end{equation}
Computing $S(x)$ at every time step is tractable: since $F(\cdot,x)$
is affine and strongly monotone and $Z$ is a polyhedron, the v-GNE can
be computed by centralized projection-type methods for strongly
monotone variational inequalities, such as projected gradient \cite{FaccPang1}. Alternatively, it can be
computed in a distributed manner, in which each player updates only
its own decision variable, by gradient-play and learning algorithms
with convergence rate guarantees for strongly monotone games
\cite{TatNedShi21}, including settings with
coupling constraints \cite{TatKamECC25}.

Since $F$ is affine in $x$ and strongly monotone in $u$, the solution
mapping is globally Lipschitz continuous
\cite[Thm.~2.3.3]{FaccPang1}. In particular,
\begin{equation}
    \|S(x)-S(y)\|
    \le
    L_S\|x-y\|,
    \qquad
    L_S\le\frac{\|G\|}{\mu_F}.
    \label{eq:S_lipschitz}
\end{equation}

%====================================================================
\section{Stability Filter Design and Guarantees}
\label{sec:stability_filter}
%====================================================================

\subsection{Closed-loop equilibrium}
%====================================================================

We first establish existence of an equilibrium of the nominal
receding-horizon closed loop. Thus, existence of such an equilibrium
does not need to be imposed separately.

\begin{proposition}%[Existence of a closed-loop equilibrium]
\label{prop:closed_loop_equilibrium}
Suppose Assumption~\ref{ass:main} holds. Then there exists
$(\bar x,\bar u)$ satisfying
\begin{equation}
    \bar u=S(\bar x),
    \qquad
    \bar x=A\bar x+B\Xi\bar u.
    \label{eq:closed_loop_equilibrium}
\end{equation}
\end{proposition}

\begin{proof}
Since $A$ is Schur, $I-A$ is nonsingular. Define
$T(u)
    =
    S\!\left((I-A)^{-1}B\Xi u\right).
    $
The mapping $T$ is continuous because $S$ is  (see~\eqref{eq:S_lipschitz}). Moreover,
$T(Z)\subseteq Z$ as any GNE solution is also feasible.

Since $Z$ is nonempty, compact, and convex, Brouwer's fixed-point
theorem implies that there exists $\bar u\in Z$ such that
$T(\bar u)=\bar u$. Define $\bar x
    =
    (I-A)^{-1}B\Xi\bar u$.
Then $\bar u=S(\bar x)$ and
$(I-A)\bar x=B\Xi\bar u$,
or equivalently
$\bar x=A\bar x+B\Xi\bar u$.
This proves \eqref{eq:closed_loop_equilibrium}.
\end{proof}
% \AB{Justify why other equilibria than $(0,0)$ may exist, I suspect they lie on the boundary of $Z$.
% Note that if $h=0$ and $Z$ does \textbf{not} contain the origin, the proposition still proves existence of a nonzero equilibrium $(\bar x, \bar u)$, which also, I suspect, lies on the boundary of $Z$.
% Claude suggest the following proposition:}
% {\color{red}
% \begin{proposition}
% Suppose Assumption~\ref{ass:main} holds with $h=0$ and $0\in Z$, and let
% $M=(I-A)^{-1}B\Xi$. Then $(0,0)$ satisfies \eqref{eq:closed_loop_equilibrium}.
% If, in addition, $u^\top(H+GM)u>0$ for all $u\in Z\setminus\{0\}$, which holds
% in particular when $\|GM\|<\mu_F$, then $(\bar x,\bar u)=(0,0)$ is the unique
% pair satisfying \eqref{eq:closed_loop_equilibrium}.
% \end{proposition}
% \begin{proof}
% $S(0)=0$ since $\langle H0, z-0\rangle=0$ for all $z\in Z$. Any pair satisfying
% \eqref{eq:closed_loop_equilibrium} has $\bar x=M\bar u$, hence
% $\langle (H+GM)\bar u, z-\bar u\rangle\geq 0$ for all $z\in Z$. Taking $z=0$
% gives $\bar u^\top(H+GM)\bar u\leq 0$, so $\bar u=0$ and $\bar x=0$.
% Finally, $\bar u^\top(H+GM)\bar u\geq(\mu_F-\|GM\|)\|\bar u\|^2$.
% \end{proof}
% }

\begin{remark}%[Multiplicity of closed-loop equilibria]
\label{rem:multiplicity}
Proposition~\ref{prop:closed_loop_equilibrium} establishes existence
of $(\bar x,\bar u)$, but Brouwer's theorem does not rule out several
fixed points of $T$. None of the results that follow, namely global
feasibility (Proposition~\ref{prop:feasibility}), nonintervention
(Proposition~\ref{prop:nonintervention}), the intervention bound
(Proposition~\ref{prop:intervention}), closed-loop stability
(Theorem~\ref{thm:stability}), or the strategic-loss analysis in
Section~\ref{sec:strategic_loss}, use uniqueness of
$(\bar x,\bar u)$. Each requires only that a single pair satisfying
\eqref{eq:closed_loop_equilibrium} to be fixed once and used consistently
as the reference in the filter constraint defined below in~\eqref{eq:filter}. When
\eqref{eq:closed_loop_equilibrium} admits multiple solutions,
($\bar x$, $\bar u$) should therefore be understood as one such pair, selected
offline (for instance via Corollary~\ref{cor:contraction_unique} below, when
applicable); every guarantee below then holds relative to that fixed
choice, uniformly in the online state $x$.
\end{remark}

\begin{corollary}[Sufficient condition for uniqueness]
\label{cor:contraction_unique}
Let
\begin{equation}
    L_T
    =
    L_S\,\|(I-A)^{-1}B\Xi\|,
    \label{eq:LT}
\end{equation}
with $L_S$ as in \eqref{eq:S_lipschitz}. If $L_T<1$, then
$(\bar x,\bar u)$ in Proposition~\ref{prop:closed_loop_equilibrium} is
the unique pair satisfying \eqref{eq:closed_loop_equilibrium}.
Moreover, for any $u^{(0)}\in Z$, the iterates
$u^{(n+1)}=T(u^{(n)})$ satisfy
\begin{equation}
    \|u^{(n)}-\bar u\|
    \le
    L_T^n\|u^{(0)}-\bar u\|,
\end{equation}
so $\bar u$, and hence $\bar x$, can be computed offline to any
desired accuracy.
\end{corollary}

\begin{proof}
For $u,v\in Z$, the Lipschitz bound \eqref{eq:S_lipschitz} gives
\begin{align*}
    \|T(u)-T(v)\|
    &=
    \Bigl\|
        S\bigl((I-A)^{-1}B\Xi u\bigr)
    \\
    &\qquad{}-
        S\bigl((I-A)^{-1}B\Xi v\bigr)
    \Bigr\|
    \\
    &\le
    L_T\|u-v\|.
\end{align*}
If $L_T<1$, $T$ is a contraction on $Z$, which is complete since it is
a compact subset of $\mathbb R^{n_u}$. The Banach fixed-point theorem
then gives a unique fixed point $\bar u\in Z$, together with the
stated geometric convergence rate of the iterates $T^n(u^{(0)})$.
Since $\bar x=(I-A)^{-1}B\Xi\bar u$, it is likewise unique.
\end{proof}

%====================================================================
\subsection{Stability filter}
%====================================================================

Since $A$ is Schur, for any $Q\succ0$ there exists a unique
$P\succ0$ satisfying the discrete-time Lyapunov equation
\begin{equation}
    A^\top PA-P=-Q.
    \label{eq:lyap}
\end{equation}
Define $V(\xi)=\xi^\top P\xi$ and
\begin{equation}
    \rho_A
    =
    \lambda_{\max}
    \left(
        P^{-1/2}A^\top PA P^{-1/2}
    \right).
    \label{eq:rhoA}
\end{equation}
By \eqref{eq:lyap},
$P^{-1/2}A^\top PAP^{-1/2}=I-P^{-1/2}QP^{-1/2}$, hence
\begin{equation}
    0\le\rho_A
    =1-\lambda_{\min}\bigl(P^{-1/2}QP^{-1/2}\bigr)<1
\end{equation}
and
\begin{equation}
    V(A\xi)\le\rho_A V(\xi)
    \qquad
    \forall \xi.
    \label{eq:A_contraction}
\end{equation}

Choose a prescribed contraction factor
\begin{equation}
    \rho\in[\rho_A,1).
    \label{eq:rho}
\end{equation}
Throughout, $\bar x$ denotes the first component of a fixed pair $(\bar x,\bar u)$ satisfying
\eqref{eq:closed_loop_equilibrium}, in the sense of
Remark~\ref{rem:multiplicity}. We define the stability-filtered
strategy by
\begin{equation}
\begin{aligned}
    \widetilde S(x)
    =
    \arg\min_{u\in Z}\quad&
        \frac12\|u-S(x)\|_R^2
        \\
    \text{s.t.}\quad&
        V(Ax+B\Xi u-\bar x)
        \le
        \rho V(x-\bar x),
\end{aligned}
\label{eq:filter}
\end{equation}
where $R\succ0$.

The filter therefore selects the feasible strategy closest to the
nominal v-GNE while enforcing a prescribed Lyapunov contraction defined by~\eqref{eq:rho}.

\begin{remark}
\label{rem:choice_Q}
The constant $\rho_A=\|P^{1/2}AP^{-1/2}\|^2$ depends on $Q$ through
$P$ and always satisfies $\rho_A\ge\rho(A)^2$, where $\rho(A)$ is the
spectral radius of $A$. Conversely, for every $\varepsilon>0$ there
exists $P\succ0$ with $\|P^{1/2}AP^{-1/2}\|\le\rho(A)+\varepsilon$,
and $Q=P-A^\top PA\succ0$ then yields
$\rho_A\le(\rho(A)+\varepsilon)^2$. Hence a suitable choice of $Q$
brings $\rho_A$ arbitrarily close to $\rho(A)^2$, which enlarges the
admissible range of contraction factors in \eqref{eq:rho}.
\end{remark}

We next establish the basic properties of the filter. We show that it
is feasible for every state (Proposition~\ref{prop:feasibility}), that
it leaves the v-GNE unchanged whenever the latter already satisfies
the contraction condition (Proposition~\ref{prop:nonintervention}),
and that its deviation from the v-GNE is bounded linearly in
$\|x-\bar x\|$ (Proposition~\ref{prop:intervention}). Building on
these results, Theorem~\ref{thm:stability} proves global exponential
stability of the filtered closed loop, and
Corollary~\ref{cor:intervention_decay} shows that the intervention
vanishes exponentially along it.

\begin{proposition}%[Global feasibility]
\label{prop:feasibility}
Under Assumption~\ref{ass:main}, the stability filter
\eqref{eq:filter} is feasible for every
$x\in\mathbb{R}^{n_x}$.
\end{proposition}

\begin{proof}
Consider the constant candidate
$u=\bar u$.
By Proposition~\ref{prop:closed_loop_equilibrium},
$\bar u=S(\bar x)\in Z$. Moreover,
\begin{align*}
    Ax+B\Xi\bar u-\bar x
    &=
    Ax+B\Xi\bar u
    -
    (A\bar x+B\Xi\bar u)
    \\
    &=
    A(x-\bar x).
\end{align*}
Consequently,
\begin{align*}
    V(Ax+B\Xi\bar u-\bar x)
    =
    V(A(x-\bar x))&\le
    \rho_A V(x-\bar x)
    \\
    &\le
    \rho V(x-\bar x).
\end{align*}
Thus $\bar u$ is feasible for \eqref{eq:filter} for every $x$.
\end{proof}

\begin{proposition}%[Nonintervention]
\label{prop:nonintervention}
If the nominal equilibrium strategy satisfies
\begin{equation}
    V(Ax+B\Xi S(x)-\bar x)
    \le
    \rho V(x-\bar x),
    \label{eq:nominal_contraction}
\end{equation}
then
$ \widetilde S(x)=S(x)$.
\end{proposition}

\begin{proof}
If \eqref{eq:nominal_contraction} holds, then $S(x)$ is feasible for
\eqref{eq:filter}. The objective in \eqref{eq:filter} is nonnegative
and equals zero at $u=S(x)$. Since $R\succ0$, this point is the unique
optimizer.
\end{proof}

\begin{proposition}%[Bound on the stability intervention]
\label{prop:intervention}
Let
\begin{equation}
    \delta(x)=\widetilde S(x)-S(x).
    \label{eq:delta}
\end{equation}
Then
$\|\delta(x)\|
    \le
    C_\delta\|x-\bar x\|$
for every $x$, where one may take
\begin{equation}
    C_\delta
    =
    \sqrt{\frac{\lambda_{\max}(R)}
             {\lambda_{\min}(R)}}\,L_S
    \le
    \sqrt{\frac{\lambda_{\max}(R)}
             {\lambda_{\min}(R)}}
    \frac{\|G\|}{\mu_F}.
    \label{eq:Cdelta}
\end{equation}
\end{proposition}

\begin{proof}
Since $\bar u$ is feasible for \eqref{eq:filter}, optimality of
$\widetilde S(x)$ gives
\[
    \|\widetilde S(x)-S(x)\|_R
    \le
    \|\bar u-S(x)\|_R.
\]
Since $\bar u=S(\bar x)$, the Lipschitz property
\eqref{eq:S_lipschitz} yields
\[
    \|\bar u-S(x)\|
    \le
    L_S\|x-\bar x\|.
\]
Using norm equivalence,
\[
    \|\delta(x)\|
    \le
    \sqrt{
        \frac{\lambda_{\max}(R)}
             {\lambda_{\min}(R)}
    }
    L_S\|x-\bar x\|,
\]
which proves the claim.
\end{proof}

\begin{theorem}%[Global exponential stability]
\label{thm:stability}
Consider the filtered closed loop
\begin{equation}
    x(t+1)
    =
    Ax(t)+B\Xi\widetilde S(x(t)),
    \label{eq:filtered_closed_loop}
\end{equation}
initialized at $x(0)=x_0$.
Then $\bar x$ is globally exponentially stable. More precisely,
\begin{equation}
       \|x(t)-\bar x\|
    \le
    \sqrt{\frac{\lambda_{\max}(P)}
             {\lambda_{\min}(P)}}
    \rho^{t/2}
    \|x_0-\bar x\|.
    \label{eq:x_decay}
\end{equation}
\end{theorem}
\begin{proof}
The contraction inequality $  V(x(t+1)-\bar x)
    \le
    \rho V(x(t)-\bar x)
$ follows directly from
the constraint in \eqref{eq:filter}.
Hence, iteration yields
$   V(x(t)-\bar x)
    \le
    \rho^t V(x_0-\bar x)$.
Finally,
\[
    \lambda_{\min}(P)\|\xi\|^2
    \le
    V(\xi)
    \le
    \lambda_{\max}(P)\|\xi\|^2
\]
gives \eqref{eq:x_decay}.
\end{proof}

\begin{corollary}%[Vanishing intervention]
\label{cor:intervention_decay}
Along the filtered closed loop,
\begin{equation}
        \|\widetilde S(x(t))-S(x(t))\|
    =
    O(\rho^{t/2}).
    \end{equation}
More precisely,
\begin{equation}
    \|\widetilde S(x(t))-S(x(t))\|
    \le
    C_\delta
    \sqrt{\frac{\lambda_{\max}(P)}
             {\lambda_{\min}(P)}}
    \rho^{t/2}
    \|x_0-\bar x\|.
\end{equation}
\end{corollary}

\begin{proof}
The result follows immediately from
Proposition~\ref{prop:intervention} and
Theorem~\ref{thm:stability}.
\end{proof}

\begin{remark}
\label{rem:closeness_vs_optimality}
Corollary~\ref{cor:intervention_decay} measures the intervention in
the space of strategies only. Closeness of $\widetilde S(x)$ to $S(x)$
does not by itself imply that $\widetilde S(x)$ is close to an
equilibrium in the strategic sense: equilibrium is characterized by
the players' unilateral optimality, i.e., by comparing
$J_i(\widetilde S(x),x)$ with the best cost player $i$ can achieve
against the filtered strategies $\widetilde S_{-i}(x)$ of the others,
rather than by the cost incurred at the v-GNE. Moreover, when
coupling constraints are active at $S(x)$, a deviation $\delta$ may
induce a regret of first order in $\|\delta\|$, whose magnitude is
governed by the associated multipliers, as opposed to the second-order
effect one would expect at an unconstrained equilibrium. This loss of
strategic optimality is quantified in the following subsection
(Theorems~\ref{thm:strategic_loss} and~\ref{thm:epsilon_decay}).
\end{remark}

%====================================================================
\subsection{Strategic loss induced by stability filtering}
\label{sec:strategic_loss}
%====================================================================

By the construction and results above, the filtered strategy $\widetilde S(x)$ is
feasible for the original game and renders the closed loop globally
exponentially stable (Theorem~\ref{thm:stability}). However, whenever
the filter intervenes, $\widetilde S(x)$ need not itself be a GNE.
What remains is therefore to characterize its optimality in terms of
the players' cost functions $J_i$, i.e., to quantify the loss of
strategic optimality incurred for the sake of stability.

For a feasible strategy $u\in Z$, define the maximum unilateral regret
\begin{equation}
    \epsilon(u,x)
    =
    \max_i
    \left\{
        J_i(u_i,u_{-i},x)
        -
        \min_{\substack{z_i:\\(z_i,u_{-i})\in Z}}
        J_i(z_i,u_{-i},x)
    \right\}.
    \label{eq:regret}
\end{equation}
Thus $u$ is an $\epsilon$-GNE whenever its maximum unilateral regret
is bounded by $\epsilon$.

For $u^\star=S(x)$, the v-GNE optimality conditions imply the
existence of common multipliers $\lambda\ge0$ and $\nu$ satisfying (see, e.g., \cite[Thm.~12.1]{FaccPang1})
\begin{subequations}
\label{eq:kkt}
\begin{align}
    F(u^\star,x)
    +C^\top\lambda
    +E^\top\nu
    &=0,
    \label{eq:kkt_stationarity}
    \\
    Cu^\star&\le d,
    \qquad
    Eu^\star=e,
    \label{eq:kkt_feasibility}
    \\
    \lambda^\top(d-Cu^\star)&=0.
    \label{eq:kkt_complementarity}
\end{align}
\end{subequations}

We first establish that the multipliers required in the subsequent
analysis can always be selected locally bounded.

\begin{lemma}
\label{lem:bounded_multipliers}
Under Assumption~\ref{ass:main}, for every compact set
$\mathcal X\subset\mathbb{R}^{n_x}$ there exists a constant
$M_\lambda(\mathcal X)<\infty$ such that, for every $x\in\mathcal X$, there exists
a KKT multiplier pair $(\lambda(x),\nu(x))$ associated with $S(x)$
satisfying
\begin{equation}
      \|\lambda(x)\|\le M_\lambda(\mathcal X).
        \label{eq:bounded_lambda}
\end{equation}
In particular, the inequality multipliers admit a locally bounded
selection in a neighborhood of $\bar x$.
\end{lemma}

The proof, which relies on Hoffman's error bound for linear
systems \cite{hoffman}, is given in
Appendix~\ref{app:proof_bounded_multipliers}.

\begin{theorem}%[Strategic-loss decomposition]
\label{thm:strategic_loss}
Let
$u^\star=S(x)$, $\tilde u=\widetilde S(x)$, and $\delta=\tilde u-u^\star$.
Under Assumption~\ref{ass:main}, let
$(\lambda(x),\nu(x))$ be any KKT multiplier pair associated with
$u^\star$. Then
\begin{align}
    \epsilon(\tilde u,x)
    &\le
    \lambda(x)^\top(d-C\tilde u)
    +
    \max_i
    \frac{\|(H\delta)_i\|^2}{2m_i} \cr
    &=  -\lambda(x)^\top C\delta
    +
    \max_i
    \frac{\|(H\delta)_i\|^2}{2m_i}.
        \label{eq:epsilon_decomposition}
\end{align}
\end{theorem}

\begin{proof}
Fix a player $i$ and any feasible unilateral deviation $z_i$, i.e.,
$z=(z_i,\tilde u_{-i})\in Z$.
By strong convexity of $J_i$ with respect to its own decision variable,
\begin{align}
    &J_i(\tilde u_i,\tilde u_{-i},x)
    -
    J_i(z_i,\tilde u_{-i},x)
    \nonumber\\
    &\qquad\le
    F_i(\tilde u,x)^\top(\tilde u_i-z_i)
    -
    \frac{m_i}{2}
    \|\tilde u_i-z_i\|^2.
    \label{eq:strong_convexity_step}
\end{align}

Since $F$ is affine, $F(\tilde u,x)
    =
    F(u^\star,x)+H\delta$.
Thus, by the v-GNE stationarity condition (see~\eqref{eq:kkt_stationarity}),
\begin{align}\label{eq:stat}
       F(\tilde u,x)
    =
    -C^\top\lambda
    -E^\top\nu
    +H\delta. 
\end{align}
Let $w=\tilde u-z$. 
Since $\tilde u$ and $z$ differ only in player $i$'s decision block,
$F_i(\tilde u,x)^\top(\tilde u_i-z_i) = F(\tilde u,x)^\top w$.
Furthermore, both $\tilde u$ and $z$ satisfy the equality constraints,
and hence
$Ew=0$.
Therefore, by~\eqref{eq:stat}
\begin{align}
    F(\tilde u,x)^\top w
    &=
    -\lambda^\top Cw
    +(H\delta)^\top w\cr
    &\le \lambda^\top(d-C\tilde u) +(H\delta)^\top w,
    \label{eq:Fw}
\end{align}
since $Cz\le d$ and $-\lambda^\top Cw
    =
    \lambda^\top(Cz-C\tilde u)$.
Substitution of~\eqref{eq:Fw} into \eqref{eq:strong_convexity_step} gives
\begin{align}\label{eq:gap1}
    &J_i(\tilde u_i,\tilde u_{-i},x)
    -
    J_i(z_i,\tilde u_{-i},x)
    \cr
    &\quad\le
    \lambda^\top(d-C\tilde u)
    +(H\delta)_i^\top(\tilde u_i-z_i)
    -
    \frac{m_i}{2}
    \|\tilde u_i-z_i\|^2\cr
    &\qquad\le
    \lambda^\top(d-C\tilde u)
    +
    \frac{\|(H\delta)_i\|^2}{2m_i}.
\end{align}
where we used $(H\delta)_i^\top(\tilde u_i-z_i)-\frac{m_i}{2}\|\tilde u_i-z_i\|^2\le \frac{\|(H\delta)_i\|^2}{2m_i}$.
Since the inequality~\eqref{eq:gap1} holds for every feasible unilateral deviation,
taking the supremum over $z_i$ and then the maximum over players gives
\[   \epsilon(\tilde u,x)
    \le
    \lambda(x)^\top(d-C\tilde u)
    +
    \max_i
    \frac{\|(H\delta)_i\|^2}{2m_i}\]

Finally,
$\lambda^\top(d-Cu^\star)=0$
implies
\begin{align*}
    \lambda^\top(d-C\tilde u)
    &=
    \lambda^\top(d-Cu^\star)
    -
    \lambda^\top C(\tilde u-u^\star)
    =
    -\lambda^\top C\delta,
\end{align*}
which proves \eqref{eq:epsilon_decomposition}.
\end{proof}

\begin{remark}
\label{rem:loss_interpretation}
The bound \eqref{eq:epsilon_decomposition} separates the strategic loss into
two qualitatively different contributions:
\begin{equation}
    \epsilon(\tilde u,x)
    \le
    \underbrace{
        \lambda(x)^\top(d-C\tilde u)
    }_{\text{constraint-release loss}}
    +
    \underbrace{
        \max_i
        \frac{\|(H\delta)_i\|^2}{2m_i}
    }_{\text{stationarity perturbation}}.
    \label{eq:loss_interpretation}
\end{equation}
The second term is quadratic in the stability intervention $\delta$.
In contrast, moving away from an inequality that is active at the
nominal v-GNE and has a positive equilibrium multiplier can generate
a first-order loss. Thus, the potentially first-order strategic
effect of stability filtering is entirely associated with releasing
positively priced active constraints.
\end{remark}

\begin{theorem}%[Global recovery of strategic optimality]
\label{thm:epsilon_decay}
Under Assumption~\ref{ass:main}, for every initial condition
$x_0\in\mathbb R^{n_x}$ there exists a constant $C_\epsilon(x_0)>0$
such that, along the filtered closed loop
\eqref{eq:filtered_closed_loop} initialized at $x(0)=x_0$,
\begin{equation}
       \epsilon(\widetilde S(x(t)),x(t))
    \le
    C_\epsilon(x_0)\|x(t)-\bar x\|,
    \qquad
    t=0,1,2,\dots
        \label{eq:epsilon_linear}
\end{equation}
Consequently,
$\epsilon(t)
    =
    \epsilon(\widetilde S(x(t)),x(t))
    =
    O(\rho^{t/2})$,
for every $x_0$. In particular, $\lim_{t\to\infty}\epsilon(t)=0$
regardless of the initial condition.
\end{theorem}

\begin{proof}
Fix $x_0$ and define the sublevel set
\begin{equation}
    \mathcal X_{x_0}
    =
    \bigl\{
        \xi\in\mathbb R^{n_x}
        \mid
        V(\xi-\bar x)\le V(x_0-\bar x)
    \bigr\}.
    \label{eq:sublevel_set}
\end{equation}
Since $V$ is a positive-definite quadratic, $\mathcal X_{x_0}$ is a
compact ellipsoid, and every $\xi\in\mathcal X_{x_0}$ satisfies
\begin{equation}
    \|\xi-\bar x\|
    \le
    R(x_0):=
    \sqrt{\frac{\lambda_{\max}(P)}
             {\lambda_{\min}(P)}}\,\|x_0-\bar x\|.
    \label{eq:Rx0}
\end{equation}
By the contraction inequality,
$V(x(t+1)-\bar x)\le\rho V(x(t)-\bar x)\le V(x(t)-\bar x)$ along the
filtered closed loop; since $x_0\in\mathcal X_{x_0}$ trivially,
induction on $t$ gives $x(t)\in\mathcal X_{x_0}$ for every
$t=0,1,2,\dots$, so $\mathcal X_{x_0}$ contains the entire trajectory.

By Lemma~\ref{lem:bounded_multipliers} applied to the compact set
$\mathcal X_{x_0}$, there is a valid KKT multiplier selection with
\[
    \|\lambda(x)\|\le M_\lambda(x_0),
    \qquad
    x\in\mathcal X_{x_0},
\]
where $M_\lambda(x_0) = M_{\lambda}(\mathcal X_{x_0})$ is a finite constant depending on $x_0$.
By Proposition~\ref{prop:intervention}, $\|\delta(x)\|\le
C_\delta\|x-\bar x\|$ for every $x$. Using
Theorem~\ref{thm:strategic_loss},
\begin{align}
    \epsilon(\widetilde S(x),x)
    &\le
    -\lambda(x)^\top C\delta(x)
    +
    \max_i
    \frac{\|(H\delta(x))_i\|^2}{2m_i}
    \nonumber\\
    &\le
    M_\lambda(x_0)\|C\|\,\|\delta(x)\|
    +
    c_H\|\delta(x)\|^2,
    \label{eq:epsilon_prebound}
\end{align}
where $c_H = \max_i\frac{\|H_i\|^2}{2m_i}$.
Therefore, for every $x\in\mathcal X_{x_0}$,
\begin{align}
    \epsilon(\widetilde S(x),x)
    &\le
    M_\lambda(x_0)\|C\|C_\delta
    \|x-\bar x\|
    +
    c_HC_\delta^2
    \|x-\bar x\|^2.
    \label{eq:epsilon_state_bound}
\end{align}
By \eqref{eq:Rx0},
$\|x-\bar x\|\le R(x_0)$ for every $x\in\mathcal X_{x_0}$. Thus, 
\[
    c_HC_\delta^2\|x-\bar x\|^2
    \le
    c_HC_\delta^2R(x_0)\|x-\bar x\|,
    \qquad
    x\in\mathcal X_{x_0}.
\]
Substituting into \eqref{eq:epsilon_state_bound} gives, for every
$x\in\mathcal X_{x_0}$ and in particular for every $x(t)$,
\[
    \epsilon(\widetilde S(x),x)
    \le
    \underbrace{
        \Bigl[
            M_\lambda(x_0)\|C\|C_\delta
            +
            c_HC_\delta^2R(x_0)
        \Bigr]
    }_{\displaystyle =:\,C_\epsilon(x_0)}
    \|x-\bar x\|,
\]
which proves \eqref{eq:epsilon_linear}.

Finally, by Theorem~\ref{thm:stability},
\[
    \|x(t)-\bar x\|
    \le
    \sqrt{\frac{\lambda_{\max}(P)}
             {\lambda_{\min}(P)}}
    \rho^{t/2}
    \|x_0-\bar x\|,
\]
so combining with \eqref{eq:epsilon_linear},
\[
    \epsilon(t)
    \le
    C_\epsilon(x_0)
    \sqrt{\frac{\lambda_{\max}(P)}
             {\lambda_{\min}(P)}}\,
    \rho^{t/2}\|x_0-\bar x\|,
\]
and therefore $\epsilon(t)\to0$, for every $x_0$.
\end{proof}

%====================================================================
\subsection{Main result}
%====================================================================

We now collect the results of this section into a single statement,
which constitutes the main result of the paper: the stability filter
renders the receding-horizon game globally exponentially stable while
preserving strategic optimality asymptotically.

\begin{theorem}%[Main result]
\label{thm:main}
Suppose Assumption~\ref{ass:main} holds. Then the following statements are true.
\begin{enumerate}
    \item \emph{Equilibrium:} there exists $(\bar x,\bar u)$ with
    $\bar u=S(\bar x)$ and $\bar x=A\bar x+B\Xi\bar u$.
    \item \emph{Global feasibility:} the filtered strategy
    $\widetilde S(x)$ defined in \eqref{eq:filter} exists for every $x\in\mathbb R^{n_x}$.
    \item \emph{Exact nonintervention:} if $S(x)$ satisfies
    \eqref{eq:nominal_contraction}, then $\widetilde S(x)=S(x)$.
    %\item \emph{Bounded intervention:}
    %$\|\widetilde S(x)-S(x)\|\le C_\delta\|x-\bar x\|$ for every
    %$x\in\mathbb R^{n_x}$, with $C_\delta$ given in \eqref{eq:Cdelta}.
    \item \emph{Global exponential stability:} for every
    $x_0\in\mathbb R^{n_x}$, the trajectory of the filtered closed
    loop $x(t+1)=Ax(t)+B\Xi\widetilde S(x(t))$, $x(0)=x_0$, satisfies
    \[
        \|x(t)-\bar x\|
        \le
        \sqrt{\tfrac{\lambda_{\max}(P)}{\lambda_{\min}(P)}}\,
        \rho^{t/2}\,\|x_0-\bar x\|,
        \qquad t\ge0.
    \]
    \item \emph{Strategic recovery:} for every
    $x_0\in\mathbb R^{n_x}$, the maximum unilateral regret along the
    filtered closed loop satisfies
    $\epsilon(t)=O(\rho^{t/2})$; in particular,
    $\lim_{t\to\infty}\epsilon(t)=0$.
\end{enumerate}
\end{theorem}
%Hence stability is enforced by construction, while the deviation from
%strategic optimality induced by the stability intervention vanishes
%exponentially as the state approaches the closed-loop equilibrium.

% {\color[rgb]{0.65,0.16,0.16}
% \section{Extensions (maybe not for ECC paper)}
% \subsection{Unstable open-loop systems}
% We extend to unstable systems, by setting $v(t)=\Xi Kx(t)+\tilde v(t)$, where $u_{\rm unc}(t)=Kx(t)$ is the Nash equilibrium of the unconstrained game. We assume $\Xi K$ is stabilizing, and $\tilde v(t)$ is the new decision variable for the game. Ideally, $\tilde v(t)\rightarrow 0$ means that the unconstrained Nash gain is recovered. The Lyapunov function is constructed for the prestabilized linear system with state-update matrix $A+B\Xi K$.

% \subsection{Set-point tracking}
% The objective is to regulate the output $y=Cx$ of the system to a desired set point $r$, rather than drive the state to the origin. To this end, we define the tracking cost function for each player as:
% \begin{equation}
%     J_i(u_i,u_{-i},x) = \!\sum_{k=0}^{N_p-1}\!{\footnotesize
%     \begin{bmatrix}u_i(k|t)-v_{r,i} \\ u_{-i}(k|t) -v_{r,-i}\\ x(k|t)-x_r\end{bmatrix}^\top\!\!\!
%     {\mathcal W}
%     \begin{bmatrix}u_i(k|t)-v_{r,i} \\ u_{-i}(k|t) -v_{r,-i}\\ x(k|t)-x_r\end{bmatrix}}
% \label{eq:tracking}
% \end{equation}
% with $x_r=Ax_r+Bv_r$, $r=Cx_r$ ($r$ = output reference). For constant $r$, this gives $h_i\neq 0$.
% %\color{black}
% }

\section{Numerical example}
\label{sec:example}
We consider a random system~\eqref{eq:system} with $n_x=2$ and $n_v=3$,
%generated by the Python Control Systems Library and scaled so that 
whose matrix $A$ has eigenvalues $0.95$ and $0.8844$. 
Each of the $N=3$ players controls one
input over the horizon $N_p=5$. After reordering the stage vector
in~\eqref{eq:Ji_quad} as $[v(k|t)^\top\ x(k|t)^\top]^\top$, player $i$'s
weight is ${\mathcal W}_i=L_iL_i^\top+e_ie_i^\top$, where
$L_i\in\R^{5\times 2}$ has random Gaussian entries and $e_i$ is the $i$th
unit vector. The low-rank term couples each player's input to the state in
a selfish way. The set $Z$ contains the shared constraints
$|v_1(k|t)+v_2(k|t)|\leq 0.4$ and $|v_1(k|t)+v_3(k|t)|\leq 0.3$ and the
local bounds $|v_j(k|t)|\leq 1$ at all prediction steps. The game is strongly
monotone with $\mu_F=0.410$ and $3.10\leq m_i\leq 5.42$. Since $h=0$ and
$0\in Z$, $(\bar x,\bar u)=(0,0)$ satisfies~\eqref{eq:closed_loop_equilibrium}.
%(also checked numerically). 
We set $Q=I$ in~\eqref{eq:lyap}, which gives
$\rho_A=0.9680$, and take $\rho=\rho_A+0.1(1-\rho_A)=0.9712$ and $R=I$. At
each $t$, $S(x(t))$ is computed with the QP-based method of~\cite{BT26}
in NashOpt~\cite{NashOpt}; when~\eqref{eq:nominal_contraction} fails,
\eqref{eq:filter} is solved as a second-order cone program by
Clarabel~\cite{GC26}.

\begin{figure}[t]
\centering
\includegraphics[width=\columnwidth]{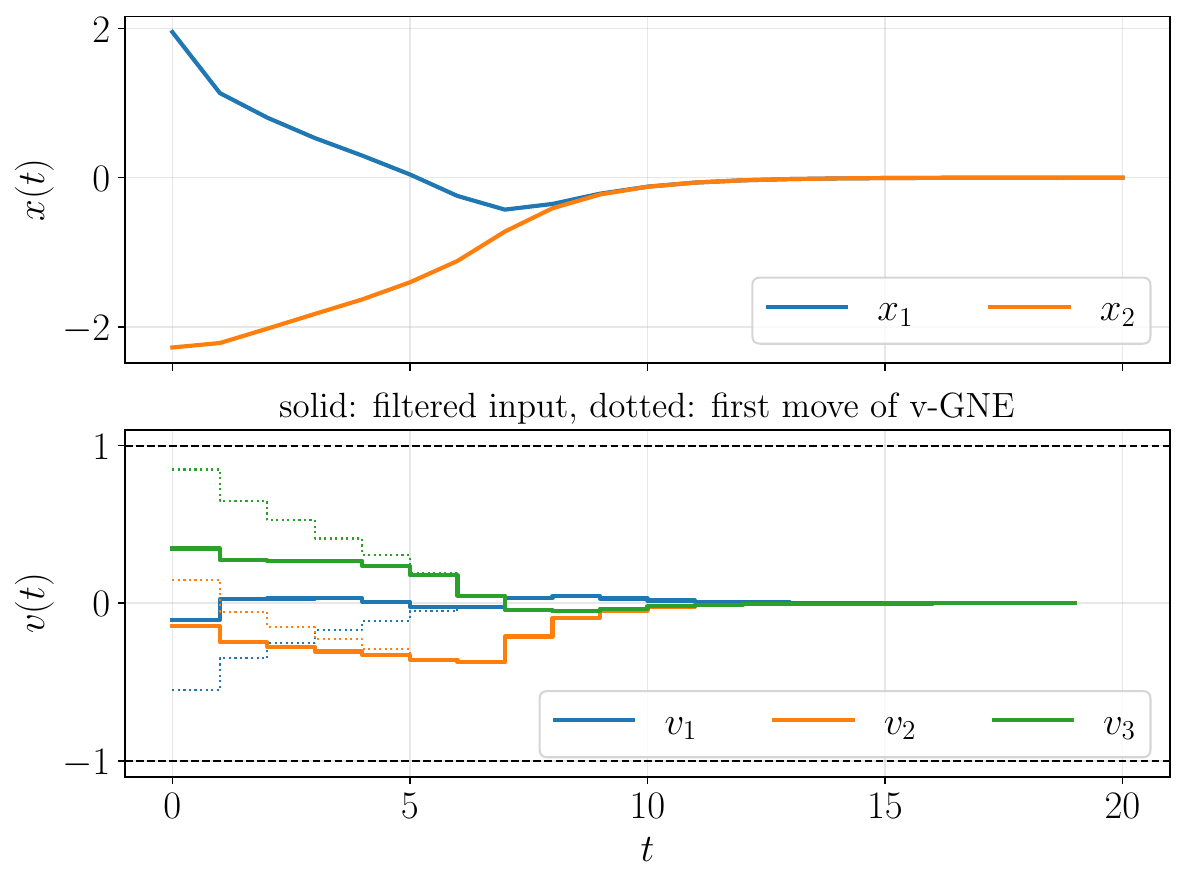}
\caption{Closed-loop states (top) and inputs (bottom). Solid: applied input
$v(t)=\Xi\widetilde S(x(t))$. Dotted: first move $\Xi S(x(t))$ of the v-GNE.
Dashed: local bounds. The red ellipse marks the steps where the filter
intervenes.}
\label{fig:states_inputs}
\end{figure}

\begin{figure}[t]
\centering
\subfigure[Phase portrait and level sets of $V$.]{%
\includegraphics[width=.48\columnwidth]{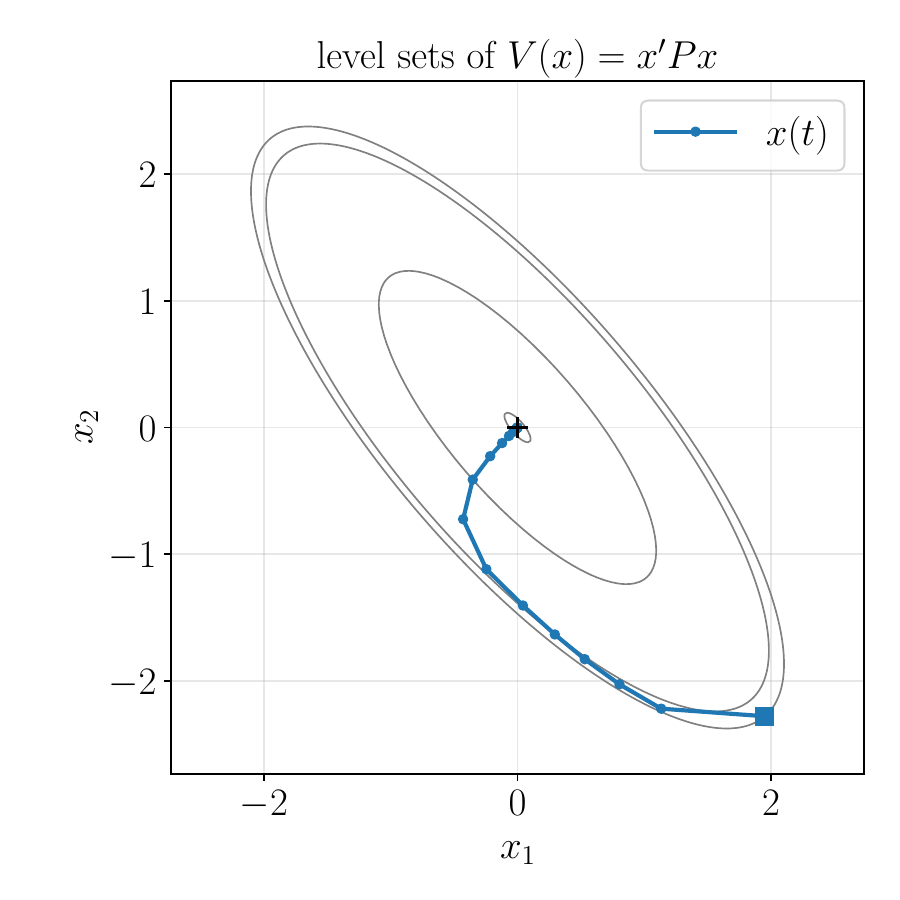}%
\label{fig:phase}}
\hfill
\subfigure[Shared constraints.]{%
\includegraphics[width=.48\columnwidth]{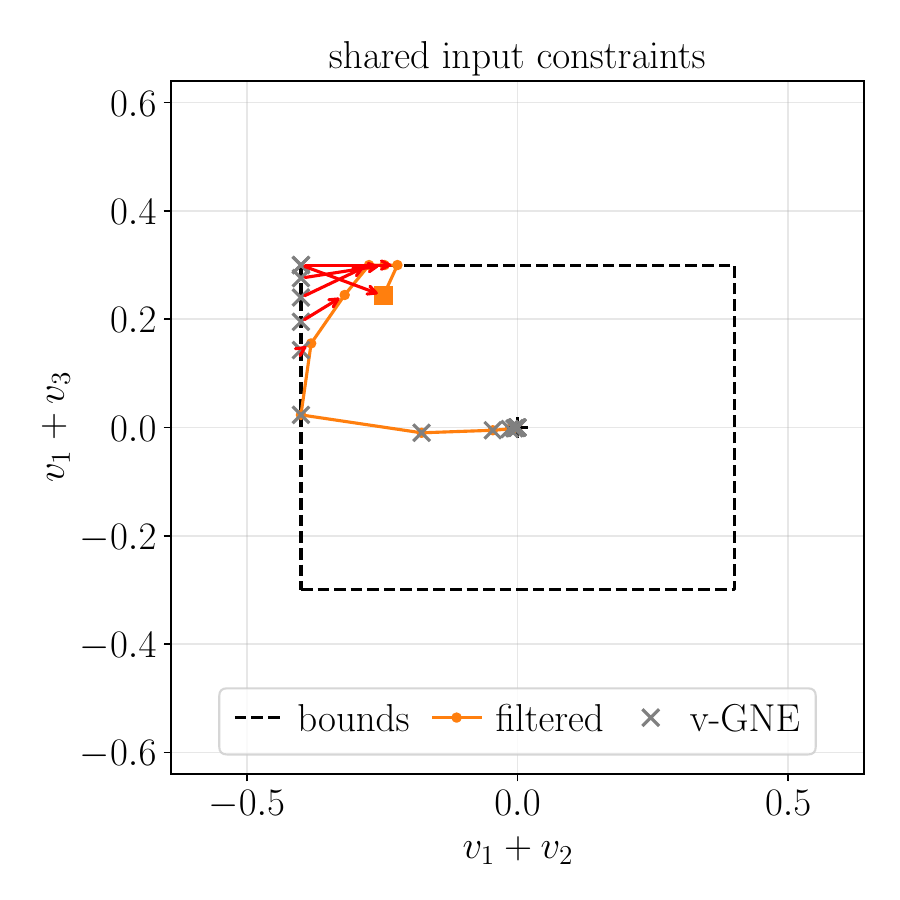}%
\label{fig:constr}}
\caption{(a) Closed-loop trajectory with level sets of $V$ through $x(t)$
every $4$ steps. (b) $v_1+v_2$ versus $v_1+v_3$ for the v-GNE first move
(crosses) and the applied input (circles), with the bounds of the shared
constraints (dashed). Red arrows: corrections applied by the filter.}
\label{fig:phase_constr}
\end{figure}

Fig.~\ref{fig:states_inputs} shows the closed loop from
$x_0=[1.9518\ \ -2.2783]^\top$. During the transient $0\leq t\leq 5$,
the v-GNE alone would increase $V$; %, with $V(x(t+1))/V(x(t))$ decreasing from $5.39$ to $1.13$. The 
the filter intervenes at these steps, with
$\|\Xi\delta(x(t))\|_\infty$ decreasing from $0.50$ to $0.03$, in agreement
with Corollary~\ref{cor:intervention_decay}. For $t\geq 6$, $S(x(t))$
satisfies~\eqref{eq:nominal_contraction} and is applied unchanged
(Proposition~\ref{prop:nonintervention}). As guaranteed by
Theorem~\ref{thm:stability}, $V(x(t+1))\leq\rho V(x(t))$ at every step, and
%$V(x(40))/V(x_0)=2.5\cdot 10^{-18}$ 
$V(x(20))\approx 0$ (Fig.~\ref{fig:phase}). During the
transient, the v-GNE lies on the edge $v_1+v_2=-0.4$ of $Z$, while the filter
moves the applied input away from it and onto the edge $v_1+v_3=0.3$ for
$1\leq t\leq 3$ (Fig.~\ref{fig:constr}).

\section{Conclusion}
\label{sec:conclusion}  
In this paper, we have proposed the use of a control-input filtering mechanism to ensure desired behavior and stability of a given constrained closed-loop game-theoretic MPC system.  We have illustrated via a simple numerical example the effectiveness of the approach in guiding the system towards the desired equilibrium while satisfying all constraints. Future work will focus on extending the methodology to control unstable open-loop systems, handling state-constraints, and solving set-point tracking problems.

%====================================================================
\appendices
\section{Proof of Lemma~\ref{lem:bounded_multipliers}}
\label{app:proof_bounded_multipliers}
%====================================================================

The proof relies on the following classical error bound for linear
systems.

\begin{lemma}[Hoffman \cite{hoffman}]
\label{lem:hoffman}
Let $A\in\mathbb{R}^{m\times n}$ and $B\in\mathbb{R}^{k\times n}$.
There exists a constant $c_{A,B}>0$, depending only on $A$ and $B$,
such that for all $a\in\mathbb{R}^m$ and $b\in\mathbb{R}^k$ for which
the polyhedron
\[
    P(a,b)=\{y\in\mathbb{R}^n\mid Ay\le a,\ By=b\}
\]
is nonempty, and for all $y_0\in\mathbb{R}^n$,
\[
    \operatorname{dist}\bigl(y_0,P(a,b)\bigr)
    \le
    c_{A,B}\bigl(\|(Ay_0-a)_+\|+\|By_0-b\|\bigr),
\]
where $(\cdot)_+$ denotes the componentwise positive part.
\end{lemma}

\begin{proofof}{Lemma~\ref{lem:bounded_multipliers}}
Fix $x\in\mathcal X$ and write $g=F(S(x),x)$ and
$\mathcal A=\{j\mid C_jS(x)=d_j\}$ for the active set.

%\emph{Step 1 (a polyhedron of multipliers).}
Consider
\[
\begin{aligned}
    \Lambda(x)=\bigl\{(\lambda,\nu)\;\big|\;
    &\lambda\ge0,\ \lambda_j=0\ (j\notin\mathcal A),
    \\
    &C^\top\lambda+E^\top\nu=-g\bigr\}.
\end{aligned}
\]
Every $(\lambda,\nu)\in\Lambda(x)$ is a KKT multiplier pair
associated with $S(x)$: \eqref{eq:kkt_stationarity} and
$\lambda\ge0$ hold by definition, \eqref{eq:kkt_feasibility} holds
since $S(x)\in Z$, and \eqref{eq:kkt_complementarity} holds since
$\lambda_j=0$ whenever $d_j-C_jS(x)>0$. Conversely, since $Z$ is
polyhedral, the KKT conditions \eqref{eq:kkt} are necessary for
$S(x)$ to solve $VI(Z,F(\cdot,x))$ \cite{FaccPang1}, so a KKT pair
exists, and by \eqref{eq:kkt_complementarity} it lies in
$\Lambda(x)$. Hence $\Lambda(x)\neq\emptyset$.

%\emph{Step 2 (applying Lemma~\ref{lem:hoffman}).}
With $y=(\lambda,\nu)$, the set $\Lambda(x)$ has the form $P(a,b)$ of
Lemma~\ref{lem:hoffman} with
\[
    A=\begin{bmatrix}-I&0\end{bmatrix},\, a=0,\,
    B_{\mathcal A}=\begin{bmatrix}C^\top&E^\top\\ I_{\mathcal A^c}&0\end{bmatrix},\,
    b=\begin{bmatrix}-g\\0\end{bmatrix},
\]
where $I_{\mathcal A^c}$ collects the rows of the identity indexed by
$j\notin\mathcal A$. The matrices $A$ and $B_{\mathcal A}$ depend only
on $C$, $E$, and $\mathcal A$; the state $x$ enters only through the
right-hand side $b$. Let $c_{\mathcal A}=c_{A,B_{\mathcal A}}$ be the
corresponding Hoffman constant.

%\emph{Step 3 (choice $y_0=0$).}
At $y_0=0$ the inequality constraints are satisfied, $(Ay_0-a)_+=0$,
and the equality residual is $\|B_{\mathcal A}y_0-b\|=\|g\|$.
Lemma~\ref{lem:hoffman} therefore gives
\[
    \operatorname{dist}\bigl(0,\Lambda(x)\bigr)\le c_{\mathcal A}\|g\|.
\]
Since $\Lambda(x)$ is nonempty and closed, it contains an element
$(\lambda(x),\nu(x))$ of minimal norm, which is a KKT multiplier pair
by Step~1 and satisfies
\[
    \|\lambda(x)\|
    \le
    \|(\lambda(x),\nu(x))\|
    =
    \operatorname{dist}\bigl(0,\Lambda(x)\bigr)
    \le
    c_{\mathcal A}\|F(S(x),x)\|.
\]

%\emph{Step 4 (uniformity in $x$).}
The active set $\mathcal A$ varies with $x$, but it is a subset of the
finitely many constraint indices. Hence
$c_Z=\max_{\mathcal A}c_{\mathcal A}<\infty$ is independent of $x$,
and
\begin{equation}
    \|\lambda(x)\|\le c_Z\|F(S(x),x)\|.
    \label{eq:lambda_bound_F}
\end{equation}
Finally, $F(S(x),x)$ is continuous with respect to $x$ (see~\eqref{eq:S_lipschitz}),
hence bounded on the compact set $\mathcal X$, which gives
\eqref{eq:bounded_lambda} with
$M_\lambda=c_Z\max_{x\in\mathcal X}\|F(S(x),x)\|$.
\end{proofof}

\bibliographystyle{IEEEtran}
\bibliography{srtrMonGames_ref}

\end{document}